\documentclass[a4paper,11pt]{article}
\usepackage{array}
\usepackage{theorem}
\usepackage{amsmath,amscd,amssymb} 
\usepackage{latexsym}
\usepackage{epsfig}

\theorembodyfont{\sl}

\newtheorem{lemma}{Lemma}[section]

\newtheorem{proposition}[lemma]{Proposition}
\newtheorem{theorem}[lemma]{Theorem}

\newtheorem{definition}[lemma]{Definition}
\newtheorem{remark}[lemma]{Remark}

\renewcommand{\AA}{\mathbb A}

\newcommand{\CC}{\mathbb C}
\newcommand{\DD}{\mathbb D}

\newcommand{\PP}{\mathbb P}

\newcommand{\TT}{\mathbb T}

\newcommand{\ZZ}{\mathbb Z}

\newcommand{\cD}{\mathcal D}
\newcommand{\cE}{\mathcal E}

\newcommand{\cO}{\mathcal O}

\newcommand{\To}{\longrightarrow}

\newcommand{\emb}{\hookrightarrow}

\renewcommand{\Tilde}{\widetilde}

\newcommand{\cross}{\times}

\newcommand{\imic}{\cong}

\newcommand{\Pic}{\mathop{\mathrm {Pic}}\nolimits}
\newcommand{\Sing}{\mathop{\mathrm {Sing}}\nolimits}
\newcommand{\Spec}{\mathop{\mathrm {Spec}}\nolimits}

\newcommand{\Aut}{\mathop{\mathrm {Aut}}\nolimits}

\newcommand{\Hom}{\mathop{\mathrm {Hom}}\nolimits}

\newcommand{\rk}{\mathop{\mathrm {rk}}\nolimits}

\newcommand{\sHom}{\mathop{\null\mathit {{\mathcal H}om}}\nolimits}

\newcommand{\bde}{\mathbf e}

\newcommand{\bdt}{\mathbf t}

\newcommand{\latt}[1]{{\langle{#1}\rangle}}
\renewcommand{\div}{\mathop{\mathrm {div}}\nolimits}

\newcommand{\qedsymbol}{\mbox{$\Box$}}
\newcommand{\qed}{\unskip\nobreak\hfil\penalty50\hskip1em\hbox{}\nobreak
\hfill\qedsymbol\parfillskip=0pt\finalhyphendemerits=0}
\newenvironment{proof}{\begin{ProofwCaption}{Proof}}{\end{ProofwCaption}}
\newenvironment{ProofwCaption}[1]
 {\addvspace\theorempreskipamount \noindent{\it #1.}\rm}
 {\qed \par \addvspace\theorempostskipamount}

\begin{document}

\title{Smooth rationally connected threefolds contain all smooth curves}
\author{G.K.~Sankaran}
\date{}
\maketitle

It is very easy to see that every smooth projective curve can be
embedded in $\PP^3$. Eisenbud and Harris asked whether the same is
true if $\PP^3$ is replaced by an arbitrary smooth rational projective
$3$-fold $X$ and Eisenbud suggested starting with the case where $X$
is toric. In that case the answer is yes, and one can see that in a
very explicit way, as was done in my preprint~\cite{Sa-pre}.

In response to~\cite{Sa-pre}, it was pointed out by J\'anos Koll\'ar
that much more is true: the property of containing every curve
sufficiently often actually characterises rationally connected
$3$-folds over the complex numbers.  In fact, this is already implicit
in the work of Koll\'ar and others on rational curves in algebraic
varieties, but had apparently not been directly noticed.

The purpose of this note is to explain these facts. In the first part
I follow Koll\'ar's hints and show how to assemble a proof of the
characterisation of rationally connected $3$-folds
(Theorem~\ref{main_RC}). In the second part, which is a shortened
version of~\cite{Sa-pre}, I show explicitly (Theorem~\ref{main_toric})
how to construct an embedding of a given curve into a given smooth
projective toric $3$-fold by toric methods.

\smallskip

{\fontsize{9pt}{11pt}\selectfont
\noindent\emph{Acknowledgements:} Much of this paper is really due to other
people. David Eisenbud asked me the question and drew
Koll\'ar's attention to my partial solution. Dan Ryder listened
patiently to me while I tried to answer the toric version. The toric
case uses ideas developed long ago in conversation with Tadao
Oda. Most importantly, J\'anos Koll\'ar pointed out in a series of
increasingly simple emails how to obtain better results, and then
allowed me to use his ideas. 
I thank all of them, and also the several people who, by asking me
about~\cite{Sa-pre}, encouraged me to write this version.
}

\section{Rationally connected varieties}

In this section $X$ is a smooth projective variety over an
algebraically closed field.

\subsection{RC and SRC}

We recall some standard definitions from~\cite{RCB} and~\cite{AK}.

\begin{definition}\label{define_SRC}\cite[IV.3.2.3]{RCB}
$X$ is \emph{separably rationally connected}, abbreviated \emph{SRC},
if there exists a variety $W$ and a morphism $u\colon \PP^1\cross
W\to X$ such that
$$
u^{(2)}\colon
(\PP^1\cross W)\cross_W(\PP^1\cross W)\imic \PP^1\cross\PP^1\cross W \To
X\cross X
$$
is dominant.
\end{definition}
In other words, $X$ is SRC if for all $x_1,\ x_2$ in some
Zariski-dense subset of $X$ we can find $w\in W$ and $P_1,\ P_2\in
\PP^1$ such that $u(P_i,w)=x_i$ for $i=1,\ 2$.

\begin{definition}\label{define_RC}\cite[IV.3.2.2]{RCB}
$X$ is \emph{rationally connected}, abbreviated \emph{RC}, if there
  exists a variety $W$, a family $U\to W$ and a morphism $u\colon U\to
  X$ such that
$$ 
u^{(2)}\colon U\cross_W U \To X\cross X
$$
is dominant.
\end{definition}
Clearly SRC$\implies$RC, and the converse is also true in
characteristic zero (\cite[IV.3.3.1]{RCB}).

\begin{definition}\label{define_very_free}\cite[Definition 8]{AK} 
A morphism $f\colon \PP^1\to X$ is said to be \emph{very free} if $f^*
T_X$ is an ample vector bundle.
\end{definition}
Recall that a vector bundle $\cE$ on $\PP^1$ is ample if and only if
$\cE=\bigoplus \cO(a_j)$ with all $a_j>0$.
\begin{lemma}\label{amplepullback}\cite[IV.3.9]{RCB}
If $X$ is a smooth projective SRC variety then there exists a very
free map $g_0\colon \PP^1\to X$.
\end{lemma}
\begin{lemma}\label{noH1}
If $X$ is a smooth projective SRC variety, then for any smooth
projective curve $C$ there is a map $g\colon C \to X$ such that
$H^1\big(g^*T_X(-P-Q)\big)=0$ for any two distinct points $P,\ Q\in C$.
\end{lemma}
\begin{proof}
Suppose $C$ has genus $g$. We choose a map $g_1\colon \PP^1\to \PP^1$
such that $g_1^*g_0^*T_X$ is sufficiently ample: it is enough to
require that $g_1^*g_0^*T_X\imic \bigoplus\cO(a_i)$ with each
$a_i>2g$, which can be achieved by taking $g_1$ to have sufficiently
high degree. Take any surjection $g_2\colon C \to \PP^1$ and let $F$
be a fibre of $g_1g_2\colon C\to \PP^1$. If we put $g=g_0g_1g_2\colon
C\to X$, we have
$$
g^*T_X(-P-Q)=\bigoplus\cO_C(a_iF-P-Q).
$$ 
This is the direct sum of line bundles of degree $>2g-2$ and therefore
$H^1\big(g^*T_X(-P-Q)\big)=0$.
\end{proof}

\begin{proposition}\label{weak_SRC}
If $C$ is any smooth projective curve and $X$ is any smooth SRC
projective variety of dimension $\ge 3$, then $C$ can be embedded in
$X$.
\end{proposition}
\begin{proof}
Choose $g\colon C\to X$ as in Lemma~\ref{noH1}. According to
\cite[II.1.8.2]{RCB} (with $B=\emptyset$), a general deformation of
$g$ is an embedding.
\end{proof}

Over the complex numbers we can show more.  

\begin{lemma}\label{very_free_locus}
Let $X$ be any smooth quasi-projective variety over $\CC$, and suppose
$x\in X$. Then there exists a subset $X_1(x)\subset X$, the complement
of a countable union of Zariski-closed sets, such that if $y\in X_1(x)$
and if the image of $f\colon \PP^1\to X$ passes through both $x$ and
$y$ then $f$ is very free.
\end{lemma}
\begin{proof}
This follows from \cite[Proposition 13]{AK}, exactly as \cite[Remark
  10]{AK} follows from~\cite[Proposition 10]{AK}. We consider one of
the countably many irreducible components $R$ of
$\Hom_x(\PP^1,X)=\{f\colon \PP^1\to X\mid f(0)=x\}$ and the evaluation
morphism $u_R\colon \PP^1\cross R\to X$ given by
$u_R(P,f)=f(P)$. The morphisms that are not very free form a closed
subscheme $R'\subseteq R$: but $u_R|_{\PP^1\cross R'}$ is not
dominant because of \cite[Proposition 13(2)]{AK}, so its image lies in
a proper closed subset $X_R\subset X$. So any $f$ that is not very
free has image contained in some $X_R$, and we take $X_1=X\setminus
\bigcup_RX_R$.
\end{proof}
This yields a characterisation of RC varieties of dimension $\ge 3$ in
terms of maps from curves.
\begin{theorem}\label{main_RC} 
 If $X$ is a smooth projective variety of dimension $\ge 3$ over
 $\CC$, then $X$ is rationally connected if and only if the following
 holds: for every smooth projective curve $C$ and zero-dimensional
 subscheme $Z\subset C$, and every embedding $f_Z\colon Z\emb X$, there
 is an embedding $f_C\colon C\emb X$ such that $f_C|_Z=f_Z$.
\end{theorem}
\begin{proof}
One direction is trivial: if every $f_Z$ extends then taking
$Z=\PP^1$ and $Z=\{0,1\}$ we recover the definition of RC.

Conversely, suppose that $X$ is RC of dimension at least~$3$ and that
$Z=\{P_1,\ldots,P_n\}$, and write $x_i=f(P_i)$. If $Z=\emptyset$ then
the statement reduces to Proposition~\ref{weak_SRC}. Otherwise, choose
$x_0\in X_1(x_1)$ as in
Lemma~\ref{very_free_locus}. By~\cite[(4.1.2.4)]{lowdeg} we can find a
map $f_0\colon \PP^1\to X$ such that $x_0,\ldots,x_n$ are all in the
image of $f_0$. See also~\cite[(5.2)]{lowdeg}. If $Z$ is not reduced,
then some of the $P_i$ are not closed, but this makes no difference
because \cite[(4.1.2.4)]{lowdeg} allows us to specify the Taylor
expansion of $f_0$ as far as we like.

The map $f_0$ is very free by Lemma~\ref{very_free_locus}. Exactly as
in Lemma~\ref{noH1} we may compose $f_0$ with suitable maps $f_2\colon
C\to \PP^1$ and $f_1\colon \PP^1\to \PP^1$ so as to get a map $f\colon
C\to X$ such that $f|_Z=f_Z$ and $\dim H^1\big(C, f^*T_X\otimes
I_Z(-P-Q)\big)=0$ for every $P,\ Q\in C$. Although $f$ need not be an
embedding, we may take $f_C$ to be a general deformation of $f$
preserving $f|_Z$, and this is an embedding by~\cite[II.1.8.2]{RCB}.
\end{proof}
\begin{remark}\label{reduced_curve}
The condition in Theorem~\ref{main_RC} that $C$ be smooth is not
strictly necessary. It is enough for $C$ to be a reduced curve whose
singularities have embedding dimension $\le \dim X$.
\end{remark}
Indeed, let $\nu\colon \Tilde C\to C$ be the normalisation. Consider
the subscheme $Z'=Z\cup \Sing C\subset C$ and let $\Tilde Z$ be the
subscheme of $\Tilde C$ given by $I_{\Tilde
  Z}=\sHom_{\cO_C}(\cO_{\Tilde C},I_{Z'})\cdot \cO_{\Tilde C}$. If
$f_{Z'}$ is an embedding of $Z'$ in $X$, extending $f_Z$, then the
argument above allows us to extend $f_{\Tilde Z}=f_{Z'}\nu\colon
{\Tilde Z}\to X$ to $f_{\Tilde C}\colon \Tilde C\to X$ in such a way
that $f_{\Tilde C}$ is an embedding away from ${\Tilde Z}$. The image
of $f_{\Tilde C}$ is then isomorphic to~$C$.

\section{Toric varieties}

In this section we look at the particular case in which $X$ is a
smooth projective toric $3$-fold over $\CC$. As toric varieties are
rational, they are in particular SRC, so by Proposition~\ref{weak_SRC}
a smooth projective toric $3$-fold contains all curves. However, in
the toric case it is possible to give a more direct proof, and one
that shows rather more concretely how to construct an embedding of a
given curve in a given toric variety $X$.

\subsection{Maps to toric varieties}

We need a good description of maps to a smooth projective toric
variety. Several descriptions available of maps to toric varieties
exist, due to Cox~\cite{Co}, Kajiwara~\cite{Ka} and others.  The
version that we use here appeared in~\cite[Section
  2]{Sa}\footnote{Warning: some other parts of~\cite{Sa} are
  incorrect.} but the proof, which is largely due to Tadao Oda, is
very short, so we give it here. We refer to~\cite{Od} for general
background on toric varieties.

Let $\Delta$ be a finite (but not necessarily complete) smooth fan for
a free $\ZZ$-module $N$ of rank $r$. Denote the corresponding toric
variety by $X$, and write $M$ for the dual lattice $\Hom(N,\ZZ)$, with
pairing $\latt{\ ,\ }\colon M\cross N\to\ZZ$. The torus is then
$\TT=\Spec(\CC[M])$, where $\CC[M]=\bigoplus_{m\in M}\CC\bde(m)$
is the semigroup ring of $M$ over $\CC$.
 
As usual, $\Delta(d)$ denotes the set of $d$-dimensional cones in
$\Delta$. For each $\rho\in\Delta(1)$, we denote by $V(\rho)$ the
corresponding irreducible Weil divisor on $X$ and by $n_\rho$ the
generator of the semigroup $N\cap\rho$. 
\begin{theorem}\label{mappingtheorem} 
  Let $Y$ be a normal algebraic variety over $\CC$. A morphism
  $f\colon Y\to X$ such that $f(Y)\cap \TT\neq\emptyset$ corresponds to
  a collection of effective Weil divisors $D(\rho)$ on $Y$ indexed by
  $\rho\in \Delta(1)$ and a group homomorphism $\varepsilon\colon M\to
  \cO_Y(Y_0)^\cross$ to the multiplicative group of invertible regular
  functions on $Y_0=Y\setminus\bigcup_{\rho\in
    \Delta(1)}D(\rho)$, such that
\begin{eqnarray}\label{intersection}
D(\rho_1)\cap D(\rho_2)\cap\cdots\cap D(\rho_s)=\emptyset\qquad
\hbox{if }\rho_1+\rho_2+\cdots+\rho_s\not\in\Delta
\end{eqnarray}
and
\begin{eqnarray}
  \label{divisor}
 \div(\varepsilon(m))=\sum_{\rho\in\Delta(1)}
\latt{m,n_\rho} D(\rho)\qquad\hbox{for all}\quad m\in M. 
\end{eqnarray}
\end{theorem}
\begin{proof}
  Suppose $f\colon Y\to X$ is a morphism with $f(Y)\cap
  \TT\neq\emptyset$.  For each $\rho\in\Delta(1)$, we take $D(\rho)$ to
  be the pull-back Weil divisor $f^{-1}(V(\rho))$, which is
  well-defined since $Y$ is normal, $X$ is smooth and $f(Y)\not\subset
  V(\rho)$. 
 
If $\rho_1+\cdots+\rho_s\not\in\Delta$, then $V(\rho_1)\cap \cdots\cap
V(\rho_s)=\emptyset$ so $D(\rho_1)\cap \cdots\cap
D(\rho_s)=\emptyset$. In this case $Y_0=f^{-1}(\TT)$ is nonempty
by assumption, and $f|_{Y_0}$ induces
\begin{eqnarray*}
f|_{Y_0}^*\colon \CC[M]\to \cO_Y(Y_0)^\cross.
\end{eqnarray*}
 The composite $\varepsilon:=f|_{Y_0}^*\circ\bde$ satisfies
\eqref{divisor}, since
\begin{eqnarray*}
\div(\bde(m))=\sum_{\rho\in\Delta(1)}\latt{m,n_\rho}V(\rho)
\qquad\hbox{for all}\quad m\in M.   
\end{eqnarray*}
 
Conversely, suppose $\{D(\rho)\}_{\rho\in\Delta(1)}$ and $\varepsilon$
satisfy \eqref{intersection} and \eqref{divisor}. For
$\sigma\in\Delta$, put
$\hat{\sigma}=\{\rho\in\Delta(1)\mid\rho\not\prec\sigma\}$. 
Then the corresponding open piece $U_\sigma$ of $X$ satisfies
\begin{eqnarray*}
U_\sigma &=&X\setminus\bigcup_{\rho\in\hat{\sigma}}V(\rho)\\
&=&\bigcap_{\rho\in\hat{\sigma}}(X\setminus V(\rho))\\
&\imic& \Spec(\CC[M\cap\sigma^{\vee}]). 
\end{eqnarray*}
Put
$Y_{\sigma}=f^{-1}(U_\sigma)=Y\setminus\bigcup_{\rho\in\hat{\sigma}}D(\rho)$. 
Then $Y=\bigcup_{\sigma\in\Delta}Y_{\sigma}$ since the $U_\sigma$
cover $X$ (or one can check this directly). 
 
For each $\sigma\in\Delta$, $M\cap\sigma^{\vee}$ is the semigroup
consisting of $m\in M$ such that $\bde(m)$ is regular on $U_{\sigma}$. 
Thus $\varepsilon(M\cap\sigma^{\vee})$ consists of regular functions
on $Y_{\sigma}$, and defines a morphism $f_{\sigma}\colon
Y_{\sigma}\to U_{\sigma}$. These morphisms glue together to give a
morphism $f\colon Y\to X$. 
\end{proof}

\subsection{Embedding a curve}

Now we apply Theorem~\ref{mappingtheorem} to the case where $Y=C$
is a smooth projective curve and $X$ is projective of
dimension~$3$. 

Let $L$ be an effective (hence ample) divisor on $C$. Let
$\Delta(1)=\{\rho_1,\ldots,\rho_r\}$, so $r-3=\rk \Pic X>0$. We write
$n_j$ and $V_j$ (rather than $n_{\rho_j}$ and $V(\rho_j)$) for the
generator and the divisor corresponding to $\rho_j\in \Delta(1)$. 

Let $\{m_1,m_2,m_3\}$ be a $\ZZ$-basis for $M$ and put
$a_{ij}=\latt{m_i,n_j}$.  The system of linear equations $\sum_{j=1}^r
a_{ij}\xi_j=0$ has rank at most $3$, so we can find nontrivial integer
solutions. In the projective case we can do better.

\begin{lemma}\label{coefficients}
  If $X$ is projective, then $\sum_{j=1}^r a_{ij}\xi_j=0$ has
  integer solutions with $\xi_j>0$ for all $j$. 
\end{lemma}
\begin{proof}
Let $H$ be a very ample divisor on $X$. We have $\sum_j a_{ij}V_j=0$
in $\Pic X$, since it is the divisor of $\bde(m_i)$. But 
$H^2V_j$ is the degree of the surface $V_j$ in the projective
embedding of $X$ under $|H|$ and is therefore positive, so it is
enough to take $\xi_j=H^2V_j$. 
\end{proof}

On $C$ we take the line bundles $\cD_j=\cO_C(\xi_j H)$, with $\xi_j$
as in Lemma~\ref{coefficients}. We may assume that $\xi_j>2g(C)$ for
all $j$, so that any nonzero linear combination of the $\cD_j$ with
nonnegative integer coefficients is very ample.

We want to specify a map $f\colon C\to X$ by means of data as in
Theorem~\ref{mappingtheorem}. Thus we must give elements $D_j$ of the
linear system $|\cD_j|$. 
\begin{lemma}\label{disjoint}
  If the $D_j$ are general in $|\cD_j|$ then they are reduced divisors
  and $\bigcap_j D_j=\emptyset$. In particular they satisfy
  \eqref{intersection} from Theorem~\ref{mappingtheorem}. 
\end{lemma}
\begin{proof}
This follows from the very ampleness of the linear systems $|\cD_j|$.
\end{proof}
To specify a map $f\colon C\to X$ we now need only choose
$\varepsilon$ according to Theorem~\ref{mappingtheorem}. This amounts
to choosing suitable trivialisations of each of the three bundles
$\cO_C(\sum a_{ij}\cD_j)$, i.e.\ non-vanishing sections of $\cO_C(\sum
a_{ij}\cD_j)$ with order $-a_{ij}$ along $D_j$. Such trivialisations
are unique up to multiplication by non-zero scalars. This means that
the map $f=f_{\DD,\bdt}$ is determined by choices of
$\DD=(D_1,\ldots,D_r)\in|\cD_1|\cross\cdots\cross|\cD_r|$ together with
a choice of an element $\bdt\in(\CC^*)^3=\TT\subset \Aut X$. In other
words, choosing the $D_j$ determines $f$ up to composition with an
element of $\TT$ acting as an automorphism of $X$. 

We note that the action of $\TT$ has no effect on the question of
whether or not the map is an embedding, and accordingly we suppress
$\bdt$ in the notation. 

Later we shall see that $f_\DD$ will turn out to be an embedding for
all sufficiently general $\DD\in|\cD_1|\cross\cdots\cross|\cD_r|$. The
next lemma shows that in order to determine whether the general
$f_\DD$ is an immersion, it is enough to check it over an affine piece
of $X$. 

\begin{lemma}\label{affinecheck}
  Suppose that, for every $\tau\in \Delta(3)$, there is a nonempty
  open subset $A_\tau \subset \prod_j|\cD_j|$ such that 
$$
f_\DD\colon C_\tau=f_\DD^{-1}(U_\tau)\to U_\tau
$$ 
is a closed immersion if $\DD\in A_\tau$. Then $f_\DD\colon C\to
X$ is a closed immersion for general $\DD\in \prod_j|\cD_j|$. 
\end{lemma}
\begin{proof}
  It is enough to take $\DD\in \bigcap_{\tau\in \Delta(3)}A_\tau$. 
\end{proof}
\begin{theorem}\label{main_toric}
  If $X$ is a projective smooth toric $3$-fold, $C$ is a smooth
  projective curve and $\cD_j$ are as above, the map
  $f_{\DD,\bdt}\colon C\to X$ is an embedding for almost all
  $\DD\in\prod_j|\cD_j|$.
\end{theorem}
\begin{proof}
  In view of Lemma~\ref{affinecheck} it remains to check that the set
  $A_\tau$ for which $f_\DD$ is an embedding above $U_\tau$ is indeed
  nonempty. 

After renumbering, we have $\tau=\rho_1+\rho_2+\rho_3$ and we consider
the semigroup $M\cap \tau^\vee$. It is generated by $l_1$, $l_2$,
$l_3\in M$ with the property that $\latt{l_i,n_i}>0$ and
$\latt{l_i,n_k}=0$ if $1\le k\le 3$ and $k\neq i$. The function
$p_i=\varepsilon_{\DD,\bdt}(l_i)=f_{\DD,\bdt}|_{C_\tau}\circ\bde(l_i)
f_{\DD,\bdt}$ is the $i$th coordinate function: it takes the value $0$
on $D_i$ and is nonzero on $D_k$ for $1\le k\le 3$, $k\neq i$. 

We first pick $D_j$ for $j>3$ once and for all, only requiring them to
be general in the sense of Lemma~\ref{disjoint}. Now choose $D_3$ so
that $D_3$ is also reduced and disjoint from the other $D_j$ chosen so
far. This is enough to determine $p_3$ up to the torus action, since
$\div(p_3)=\latt{l_3,n_3}D_3+\sum_{j>3}\latt{l_3,n_j}D_j$ is
independent of $D_1$ and $D_2$. Similarly a choice of $D_1$ or of
$D_2$ determines $p_1$ or $p_2$ up to the torus action, independently
of the choice. 

After making such a choice of $D_3$, we claim that for general $D_2\in
|\cD_2|$ the map $(p_2,p_3)\colon C_\tau\to \AA^2$ is generically
injective. We shall check this by exhibiting a choice of $D_2$ which
makes this map injective near $D_3$. Observe that for any pair $P,\ Q\in
D_3$ (so $p_3(P)=p_3(Q)=0$) we can find $D_2\in|\cD_2|$ such that
$P\in D_2$ but $Q\not\in D_2$ (although such a choice of $D_2$ will
not be general in the sense of Lemma~\ref{disjoint}), because $\cD_2$
is sufficiently ample. For this choice of $D_2$, we have $0=p_2(P)\neq
p_2(Q)$, so $p_2(P)\neq p_2(Q)$ for general $D_2$ and hence for
general $D_2$ the values of $p_2$ on the points of $D_3$ are all
different from one another. In particular $(p_2,p_3)$ corresponding to
a general $D_2$ is injective at any point of $D_3$ and is therefore
injective generically. 

By exactly the same argument, a general choice of $D_1\in |\cD_1|$
separates points not separated by the other choices. If $P'$ and $Q'$
are (possibly infinitely near) points such that $p_2(P')=p_2(Q')$ and
$p_3(P')=p_3(Q')$, then $p_1(P')\neq p_1(Q')$ if $P'\in D_1$ and
$Q'\not\in D_1$. Such $D_1$ exist if $\cD_1$ is sufficiently ample. So
for general $D_1$ we also have $p_1(P')\neq p_1(Q')$, as required. 
\end{proof}

\bibliographystyle{alpha}

\begin{thebibliography}{Ko1}
\bibitem[AK]{AK} C. Araujo \& J. Koll\'ar, {\it Rational curves on
  varieties.} Higher dimensional varieties and rational points
  (Budapest, 2001), 13--68, Bolyai Soc. Math. Stud., {\bf 12},
  Springer, Berlin, 2003.
\bibitem[Co]{Co} D.A. Cox, {\it The functor of a smooth toric
  variety.}  Tohoku Math.  J. (2) {\bf 47} (1995), 251--262.
\bibitem[Ka]{Ka} T. Kajiwara, {\it The functor of a toric variety with
  enough invariant effective Cartier divisors.} Tohoku Math. J. (2)
  {\bf 50} (1998), 139--157.
\bibitem[Ko1]{lowdeg} J. Koll\'ar, {\it Low degree polynomial
  equations: arithmetic, geometry and topology.} European Congress of
  Mathematics, Vol. I (Budapest, 1996), 255--288, Progr. Math., {\bf
    168}, Birkh\"auser, Basel, 1998.
\bibitem[Ko2]{RCB} J. Koll\'ar, {\it Rational curves on algebraic
  varieties.}  Ergebnisse der Mathematik und ihrer Grenzgebiete {\bf
  32}, Springer, Berlin, 1996.
\bibitem[Od]{Od} T. Oda, {\it Convex bodies and algebraic geometry.}
  Ergebnisse der Mathematik und ihrer Grenzgebiete {\bf 15}, Springer,
  Berlin, 1988.
\bibitem[Sa1]{Sa} G.K. Sankaran, {\it Abelian surfaces in toric
  $4$-folds.}  Math. Ann. {\bf 313} (1999), 409--427.
\bibitem[Sa2]{Sa-pre} G.K. Sankaran, {\it Any smooth toric threefold
  contains all curves.} Preprint {\tt math.AG/0710.3290}, 2007.
\end{thebibliography}

\bigskip

\noindent
G.K.~Sankaran, Department of Mathematical Sciences, University of
Bath, Bath BA2 7AY, England\\
{\tt gks@maths.bath.ac.uk}
\end{document}